\documentclass[11pt]{article}
\usepackage{a4wide}
\usepackage{amsmath}

\usepackage{amsfonts}
\usepackage{amssymb}
\usepackage{euscript}
\newenvironment{proof}{\noindent {\it Proof.~~}\ }{\  \rule{1mm}{2mm}\medskip}

\newtheorem{theorem}{Theorem}
\newtheorem{lemma}[theorem]{Lemma}
\newtheorem{corollary}[theorem]{Corollary}
\newtheorem{proposition}[theorem]{Proposition}

\newtheorem{theirtheorem}{Theorem}

\newtheorem{theirlemma}[theirtheorem]{Lemma}

\newcommand{\subgp}[1]{\langle{#1}\rangle}

\newcommand{\qed}{\hfill \mbox{$\Box$}}
\begin{document}
\title{
  Hyper-atoms and  the critical pair Theory}

\author{ Yahya O. Hamidoune\thanks{UPMC Univ Paris 06,
 E. Combinatoire, Case 189, 4 Place Jussieu,
75005 Paris, France,     {\tt hamidoune@math.jussieu.fr} }
}
\maketitle

\begin{abstract}

We introduce the notion of a hyper-atom.
One of the main results of this paper is the $\frac{2|G|}3$--Theorem:

Let $S$ be  a finite  generating subset  of
 an abelian group $G$ of order $\ge 2$.
Let $T$ be  a finite  subset  of
 $G$ such that
 $2\le |S|\le |T|$,  $S+T$ is aperiodic,  $0\in S\cap T$ and
 $$ \frac{2|G|+2}3\ge |S+T|= |S|+|T|-1.$$
 Let $H$ be a  hyper-atom of $S$. Then $S$ and $T$ are $H$--quasi-periodic.
  Moreover
 $\phi(S)$ and $\phi(T)$ are arithmetic progressions with the same
 difference, where $\phi :G\mapsto G/H$ denotes
the canonical morphism.

This result implies easily the traditional critical pair Theory and its basic
stone: Kemperman's Structure Theorem.
\end{abstract}

\section{Introduction}

For a subset $A$  of an abelian group $G$, the {\em period} of $A$ is
$\Pi (A)=\{x\in G :A+x=A\}$. The set $A$ is said to be  {\em periodic} if $\Pi(A)\neq \{0\}$.
A basic tool in Additive Number Theory is the following
generalization of the Cauchy-Davenport Theorem 
due to Kneser:

\begin{theirtheorem}[Kneser \cite{tv}]\label{kneser}
Let $G$ be an abelian group and let $A, B\subset G$ be finite
subsets of $G$ such that $A+B$ is aperiodic.  Then $|A+B|\ge |A|+|B|-1$.

\end{theirtheorem}

The description of the subsets $A$ and $B$ with $|A+B|= |A|+|B|-1$, obtained by kemperman in \cite{kempacta}  is a deep result in the classical critical pair Theory. Another step in this direction is proposed by Grynkiewicz in \cite{davkem}. The cumulated proofs of these two results is about 80 pages.
One of our  aims in the present work is to present a methodology leading to new easier results and shortest proofs for the existing ones.
This  work is essentially self-contained. We assume only Kneser's Theorem, Proposition
\ref{Cay}, Theorem \ref{2atomejc} and  Proposition \ref{strongip}. The last three results  are proved in around  3 pages  in \cite{hiso2007}.

The isoperimetric method is a global approach introduced by the author,
which
derive additive inequalities from global properties of the fragments
and atoms.
 The reader may refer to the recent paper \cite{hiso2007} for an introduction to the applications of this method.

 For a subset $X$, we put $\partial _S(X)=(X+S)\setminus X$ and $X^S=G\setminus (X+S)$.

 Suppose that $|G|\ge 2k-1$ and let  $0\in S$ be a generating subset.
 The {\em $kth$--connectivity}
of $S$
 is defined  as
$$
\kappa _k (S )=\min  \{|X+S|-|X|\   :  \ \
\infty >|X|\geq k \ {\rm and}\ |X+S|\le |G|-k\},
$$
where $\min \emptyset =|G|-2k+1$.

 We shall say that a subset
$X$ induces a {\em $k$--separation} if $ |X|\geq k$ and $|X^S|\geq
k$. We shall say that $S$ is $k$--separable if some $X$ induces a
$k$--separation.

 A finite subset $X$ of $G$ such that $|X|\ge k$,
$|X^S|\ge k$ and $|\partial (X)|=\kappa _k(S)$ is
called a {\em $k$--fragment} of $S$. A $k$--fragment with minimum
cardinality is called a {\em $k$--atom}.



Let $0\in S$ be a generating subset of an abelian group $G.$ We
shall say that $S$ is a {\em Vosper subset} if for all $X\subset G$
with $|X|\ge 2$, we have $|X+S|\ge \min (|G|-1,|X|+|S|)$.

A subgroup with maximal cardinality
which is a $1$--fragment will be called a {\em hyper-atom}.
In Section 3, we  prove the existence of hyper-atoms and obtain the following
result:

Let $S$ be a finite generating subset of an abelian group $G$ such
that $0 \in S,$  $| S | \leq (|G|+1)/2$ and $\kappa _2 (S)\le
|S|-1.$ Let $H$ be a hyper-atom of $S$. Then $\phi (S)$ is either an arithmetic progression or a Vosper
subset, where $\phi$ is the canonical morphism from $G$ onto $G/H$.

A set  $A$ is said to be {\em $H$--quasi-periodic} if there is an $x$ such that
$(A\setminus (x+H))+H=A\setminus (x+H)$. 

In Section 5, we   apply the
global isoperimetric methodology introduced in \cite{hiso2007}
to prove the following Vosper's type result:

Let $S$ be  a finite  generating subset  of
 an abelian group $G$ of order $\ge 2$.
Let $T$ be  a finite  subset  of
 $G$ such that
 $|S|\le |T|$,  $S+T$ is aperiodic,  $0\in S\cap T$ and
 $$ \frac{2|G|+2}3\ge |S+T|= |S|+|T|-1.$$
 Let $H$ be a  hyper-atom of $S$. Then $S$ and $T$ are $H$--quasi-periodic.
  Moreover
 $\phi(S)$ and $\phi(T)$ are arithmetic progressions with the same
 difference, where $\phi :G\mapsto G/H$ denotes
the canonical morphism.

 This $\frac{2|G|}3$--Theorem implies easily several critical pair results. As an illustration
 we deduce from  it  new proofs of
 Kemperman's
Structure Theorem and Lev's Theorem.

Quite likely, the methods introduced
in the present work  lead to
descriptions for subsets $A,B$ with $|A+B|=|A|+|B|+m$, with
some small other values of $m\ge 0$.  However we shall
limit ourselves to  the case $m=-1$ in order to illustrate the
method in a relatively simple context.

\section{Terminology and preliminaries}

Let  $ A$ and $B$ be subsets of $ G $. The subgroup generated by  $A$ will
be denoted by $\subgp{A}$. The {\em Minkowski sum} is defined as
$$A+B=\{x+y \ : \ x\in A\  \mbox{and}\ y\in
  B\}.$$



Recall the following two results:
\begin{theirlemma}(folklore)
Let $G$ be a finite group and let
 $A$ and $B$  be  subsets
 such that $|A|+|B|\ge |G|+1$.
 Then  $A+B=G$.

\label{prehistorical}
 \end{theirlemma}

\begin{theirtheorem}\label{scherk}(Scherk)\cite{scherck}
  Let $X$ and $Y$ be nonempty finite subsets of an abelian group $G$. If there is
an element   $c$   of  $G$ such
  that $|X\cap(c-Y)| = 1$, then
    $|X + Y| \geq |X| + |Y| - 1.$
\end{theirtheorem}

Scherck's Theorem follows easily from Kneser's Theorem.






By a {\em proper} subgroup of $G$ we shall mean a subgroup of $G$ distinct from $G$.

The next lemma is related to a notion introduced by Lee \cite{lee}:

\begin{theirlemma}\cite{balart}{Let  $X$ be a
 subset of $G$. Then $(X^S)^{-S}+S=X+S$. \label{lee}}
\end{theirlemma}

Clearly $X\subset (X^S)^{-S}$.
Take $x\notin X+S$. Then $x\notin X^S$  and hence $x-S\subset X^S-S.$ It follows that
 $x\in (X^S-S)+S=G\setminus ((X^S)^{-S}).$\qed



Throughout all this section, $S$ denotes a finite generating subset
of an abelian group $G$ with $0\in S$. Note that the best one can get, using the isoperimetric method, for a general subset $S$ is
obtained  by decomposing modulo the subgroup generated
by a translated copy of $S$ containing $0$.




 The reader may
find all basic facts from the isoperimetric method in the recent
paper \cite{hiso2007}.


Notice that
    $\kappa _k (S)$ is the maximal integer $j$
such that for every finite subset $X\subset G$  with $|X|\geq k$,

\begin{equation}
|X+S|\geq \min \Big(|G|-k+1,|X|+j\Big).
\label{eqisoper0}
\end{equation}

Formulae (\ref{eqisoper0}) is an immediate consequence of the
definitions. We shall call (\ref{eqisoper0}) the {\em isoperimetric
inequality}. The reader may use the conclusion of this lemma as a
definition of $\kappa _k (S)$. Since $|\partial (\{0\})|\ge \kappa _1$, we have
\begin{equation}\label{bound}
\kappa _1(S)\le |S|-1.
\end{equation}

The basic intersection theorem is the following:

\begin{theorem}\cite{halgebra,hiso2007}
Let $S$ be generating subset of an abelian group $G$ with $0\in S$.
 Let $A$ be a $k$--atom and let
   $F$   be a   $k$-fragment such that  $|A\cap F|\ge k$. Then
  $A\subset F.$
In particular,  distinct $k$-atoms intersect in  at most $k-1$
elements.

\label{inter2frag} \end{theorem}

The structure of $1$--atoms is the following:

\begin{proposition} \label{Cay}\cite{hejc2,hjct}

Let  $ S$  be a  generating subset of an abelian  group $G$ with
$0\in S$.  Let $H$ be a $1$--atom of $S$ with $0\in H$.
   Then
   $H$ is a subgroup. 
 Moreover
 \begin{equation}\label{olson}
\kappa _1(S)\geq \frac{|S|}{2}.
\end{equation}

\end{proposition}

\begin{proof}
Take $x\in H$. Since $x\in (H+x)\cap H$ and since $H+x$ is a
$1$--atom, we have $H+x=H$ by Theorem \ref{inter2frag}. Therefore
$H$ is a subgroup. Since $S$ generates $G$, we have $|H+S|\ge 2|H|$,
and hence  $\kappa _1(S)=|H+S|-|H|\ge \frac{|S+H|}{2}\ge
\frac{|S|}{2}.$
\end{proof}

Recently, Balandraud introduced some isoperimetric objects and
proved  a strong form of Kneser's Theorem using {Proposition} \ref{Cay}.




 The next result is proved in \cite{Hejcvosp1}. The finite case is reported with almost the same proof in
 \cite{hactaa}.

\begin{theorem} { \cite{Hejcvosp1,hactaa}\label{2atom} Let  $S$  be a finite generating
 $2$--separable subset of an abelian group $G$ with $0\in S$ and $\kappa _2 (S) \leq |S|-1$.
Let
 $H$ be a $2$--atom with $0\in H$. Then $H$ is  a subgroup  or $|H|=2.$

\label{2atomejc} }
\end{theorem}
A short proof of this result is given in \cite{hiso2007}.

\begin{corollary} { [\cite{Hejcvosp1},Theorem 4.6]\label{ejcf}
Let  $S$  be a $2$--separable finite
subset of an
 abelian group $G$ such that $0\in S$, $|S|\leq (|G|+1)/2$ and $\kappa _2 (S) \leq |S|-1$.

 If  $S$ is not an arithmetic progression then there
 is a subgroup $H$ which is a $2$--fragment of $S$.

\label{vosper}                      }
\end{corollary}

\begin{proof}

Suppose that $S$ is not an arithmetic progression.

 Let $H$ be a  $2$--atom such that $0\in H$. If $\kappa
_2\leq |S|-2$, then clearly $\kappa _2=\kappa _1$ and $H$
is also a $1$--atom. By Proposition \ref{Cay}, $H$ is a subgroup.
Then we may assume $$\kappa _2(S)=|S|-1.$$  By Theorem
\ref{2atomejc}, it would be enough to consider the case $|H|=2$, say
$H=\{0,x\}$. Put $N=\subgp{x}.$

Decompose $S=S_0\cup \cdots \cup S_j$ modulo $N$, where $|S_0+H|\le
|S_1+H| \le \cdots \le |S_j+H|.$ We have $|S|+1=|S+H|=\sum
\limits_{0\le i \le j}|S_i+\{0,x\}|.$

Then $|S_i|=|N|$, for all $i\ge 1$.  We have $j\ge 1$, since
otherwise $S$ would be an arithmetic progression. In particular $N$
is finite.
 We have
$|N+S|<|G|$, since otherwise   $|S|\ge |G|-|N|+1\ge
\frac{|G|+2}{2},$ a contradiction.

Now  \begin{eqnarray*} |N|+|S|-1&=&|N|+\kappa _2(S)\\&\le& |S+N|=
(j+1)|N|\\&\le& |S|+|N|-1, \end{eqnarray*}

 and hence $N$ is a $2$-fragment.
\end{proof}

 Corollary \ref{vosper} was used to solve Lewin's
Conjecture on the  Frobenius number
 \cite{hactaa}.
Corollary \ref{vosper} coincides with [\cite{Hejcvosp1},Theorem
4.6]. A special case of this result is Theorem 6.6 of \cite{hactaa}.
As mentioned in \cite{hplagne}, there was a misprint in this last
statement. Indeed $|H| + |B| - 1$ should be replaced by $|H| + |B|$
in case (iii) of [ Theorem 6.6, \cite{hactaa}].

Alternative proofs of Corollary \ref{vosper} (with $|S|\leq |G|/2$
replacing $|S|\leq (|G|+1)/2$),  using Kermperman's Structure Theorem, were
obtained by Grynkiewicz in \cite{davdecomp} and Lev in
\cite{levkemp}. In the present paper, Corollary \ref{vosper} will be
one of the pieces leading to a new proof of Kemperman's Theorem.


Let  $H$ be a subgroup. A partition  $A=\bigcup \limits_{i\in I}
A_i$  will be called
 a $H$--{\em decomposition}
of $A$ if for every $i$, where $A_i$ is the nonempty intersection of some $H$--coset
with $A$.
A $H$--decomposition  $A=\bigcup \limits_{i\in I}
A_i$  will be called
 a $H$--{\em modular-progression } if it is an arithmetic progression modulo $H$.

  We need the following consequence of
Menger's Theorem:

\begin{proposition} \cite{hiso2007}{ Let $G $ be an abelian group and let $S$ be a finite subset of $G$ with $0\in S$.
Let $H$ be a subgroup of $G$
 and let   $S=S_0\cup \cdots \cup S_u$  be a $H$-decomposition with $0\in S_0$.
 Let   $X=X_0\cup \cdots \cup X_t$  be a $H$-decomposition with $\frac{|G|}{|H|}\ge t+u+1.$
 Assume that $\kappa _1(\phi (S))\ge u.$

 Then there are
pairwise distinct elements
 $n_1, n_2, \cdots, n_{r} \in [0,t]$ and  elements
 $y_1, y_2, \cdots, y_{r} \in S\setminus H$ such that

 $$|\phi(X\cup X_{{n_1}}y_1\cup \cdots  \cup X_{{n_r}}y_{r})|=t+u+1.$$

\label{strongip}}
\end{proposition}

We call the property given in Proposition  \ref{strongip}  the {\em
strong isoperimetric property}.

\section{Hyper-atoms }

In this section, we investigate  the new notion of a hyper-atom.
Recall that $S$ is a Vosper subset if and only if  $S$ is non
$2$--separable or if $\kappa _2(S)\ge |S|$.

\begin{lemma} { Let $S$ be a finite generating Vosper subset of an
abelian group $G$ such that $0 \in S$.   Let $X\subset G$  be such
that   $|X+S|=|X|+|S|-1$. Also assume that $|X|\ge |S|$ if $|X|+|S|=|G|$. Then for every $y\in
S$, we have $|X+(S\setminus \{y\})|\ge |X|+|S|-2$. \label{vominus}}
\end{lemma}
\begin{proof}

By  the definition of a Vosper subset, we have $|X+S|\ge |G|-1$.
There are   two possibilities:

{\bf Case} 1. $|X+S|= |G|-1$.

 Suppose that $|X+ (S\setminus \{y\})|\le |X|+|S|-3$
and take an element  $z$ of $(X+S)\setminus (X+(S\setminus \{y\}))$.
We have $z-y\in X$. Also $(X\setminus \{z-y\})+S\subset
((X+S)\setminus \{z\})$. By the definition of
a Vosper subset, we have $|(X\setminus \{z-y\})+S|\ge \min
(|G|-1,|X|-1+|S|)=|X|+|S|-1$. Clearly $X+S\supset ((X\setminus
\{z-y\})+S)\cup \{z\}$. Hence $|X+S|\ge |X|+|S|$, a contradiction.

{\bf Case} 2. $|X+S|= |G|$.

 Suppose that $|X+ (S\setminus \{y\})|\le |X|+|S|-3$
and take a $2$--subset   $R$ of $(X+S)\setminus ( X+(S\setminus
\{y\}))$. We have $R-y\subset X$. Also $(X\setminus (R-y))+S\subset
(X+S)\setminus R$. By the definition of a
Vosper subset, $|(X\setminus (R-y))+S|\ge \min (|G|-1,|X|-2+|S|)$.
We have  $|X|=1$. Otherwise and since  $X+S\supset ((X\setminus
(R-y))+S)\cup R$, we have  $|X+S|\ge |X|+|S|$,  a contradiction.
Then $|X|=|S|=3$, and hence $|G|=5$. Now
by the Cauchy Davenport Theorem, $|X+(S\setminus \{y\})|\ge
|X|+|S|-2$, a contradiction.
\end{proof}

Let us prove a lemma about the {fragments in quotient groups}.

\begin{lemma} { Let $G $ be an abelian group and let $S$ be a
    finite  $2$-separable generating subset containing $0$. Let $H$ be a
subgroup which is a $2$--fragment and let $\phi : G\mapsto G/H$ be
the canonical morphism. Then
\begin{equation}\label{cosetgraph}
\kappa _1(\phi (S))=  |\phi (S)|-1.
\end{equation}
 Let $K$ be a subgroup which is a $1$--fragment of
$\phi (S)$. Then $\phi ^{-1}(K)$ is a $2$--fragment of $S$.
}\label{quotient}
\end{lemma}
\begin{proof}

 Put $|\phi (S)|=u+1$.  Since $|G|>|H+S|,$ we have $\phi (S)\ne G/H$,
 and hence $\phi (S)$ is
$1$--separable.

Let $X\subset G/H$ be such that  $X+\phi (S)\neq G/H$. Clearly
$\phi^{-1} (X)+S\neq G$. Then $|\phi^{-1} (X)+S|\ge |\phi^{-1}
(X)|+\kappa _2(S)= |\phi^{-1} (X)|+u|H|.$

It follows that $|X+\phi (S)||H|\ge |X||H|+u|H|.$ Hence $\kappa
_1(\phi (S))\ge u=|\phi (S)|-1$. The reverse inequality is obvious
and follows by (\ref{bound}). This proves (\ref{cosetgraph}).

Let $K$ be a subgroup which is a $1$--fragment of $\phi (S)$. Then
$|K+\phi (S)|=|K|+u$. Thus $|\phi ^{-1} (K)+S|=|K||H|+u|H|.$ In
particular, $\phi ^{-1} (K)$ is a $2$--fragment.\end{proof}


Let $S$ be a finite generating subset of an abelian group $G$ such
that $0 \in S.$ Proposition \ref{Cay} states that there is a $1$--atom
of $S$ which is a subgroup. A subgroup with maximal cardinality
which is a $1$--fragment will be called a {\em hyper-atom} of $S$. This
definition may be adapted to non-abelian groups and even abstract
graphs. As we shall see, the hyper-atom is more closely related to
the critical pair theory than the $2$--atom.

\begin{theorem}\label{hyperatom}
Let $S$ be a finite generating subset of an abelian group $G$ such
that $0 \in S,$  $| S | \leq (|G|+1)/2$ and $\kappa _2 (S)\le
|S|-1.$ Let $H$ be a hyper-atom of $S$. Then

$\phi (S)$ is either an arithmetic progression or a Vosper
subset, where $\phi$ is the canonical morphism from $G$ onto
$G/H$.
\end{theorem}

\begin{proof}

Let us show that \begin{equation}\label{referee}
2|\phi (S)|-1\le \frac{|G|}{|H|}.\end{equation}
 Clearly we may
assume that $G$ is finite.

Observe that $2|S+H|-2|H|\le 2|S|-2< |G|.$ It follows, since$|S+H|$
is a multiple of $|H|$, that $2|S+H|\le  |G|+|H|,$ and hence (\ref{referee}) holds.

 Suppose now that  $\phi (S)$ is not a Vosper subset. By the
 definition of a Vosper subset,  $\phi (S)$ is $2$--separable and $\kappa _2(\phi(S))\le |\phi(S)|-1.$

 Observe that  $\phi(S)$ can not
have  a $1$--fragment $M$ which is a non-zero subgroup. Otherwise  by Lemma
\ref{quotient}, $\phi ^{-1}(M)$ is a $2$--fragment of
$S$ containing strictly $H$, contradicting the maximality of $H$. By (\ref{referee}) and
Corollary \ref{vosper}, $\phi(S)$ is an arithmetic progression.\end{proof}

Theorem \ref{hyperatom} implies  a result proved by Plagne and the author \cite{hplagne}
and some extensions of it,  proved using
Kermperman's Theory,  obtained  by Grynkiewicz in
\cite{davdecomp} and Lev in \cite{levkemp}.

All these results follow from Theorem \ref{hyperatom}. Two main new facts
in Theorem \ref{hyperatom} are:
\begin{itemize}
\item
The subgroup $H$ in Theorem \ref{hyperatom}  is well described
as a hyper-atom;
\item
The equality $|H+S|-|H|=\kappa _1$ is   much precise
than the inequality $|H+S|\le |H|+|S|-1$ in the previous results.  This equality will be
needed later.
\end{itemize}


\section{Transfer Lemmas}

\begin{lemma}\label{AP1}
Let $G$ be an abelian group and
let $Y\subset G$ be an arithmetic progression containing $0$. Let $X\subset \subgp{Y}$
such that $|X+Y|=|X|+|Y|-1$.
Suppose that $|\subgp{Y}|\neq |X+Y|$ or $|Y|=2$

Then
$X$ is an arithmetic progression with the same difference as $Y$.
\end{lemma}

\begin{proof}
 The Lemma is obvious, once observed that
 a subset with cardinality $|\subgp{Y}|-1$ of the cyclic group
 $\subgp{Y}$ (generated by an arithmetic progression containing $0$)
 is an arithmetic progression  with arbitrary difference.\end{proof}


\begin{lemma}\label{nongenerating}
Let $S$ and $T$ be finite subsets of an abelian group $G$ such that
 $0\in S\cap T$ and $S+T$ is aperiodic.
  Also assume that  $2\le |S|$ and that $|S+T|= |S|+|T|-1$.
  If  $T\not\subset \subgp{S}$
then $T$ is $M$--quasi-periodic, where   $\subgp{S}=M$.
\end{lemma}

\begin{proof}
Decompose $T=T_1\cup \cdots \cup T_t$ modulo $M$.  Put $W=\{i \in [0,t] :  |T_i+S|<|M|\}.$
Therefore we have using (\ref{olson}),

\begin{eqnarray*}|T+S|&=& \sum _{i\in W}|T_i+S|+\sum _{i\notin W}|M|\\&\ge&
\sum _{i\in W}(|T_i|+\kappa _1(S))+ \sum _{i\notin W}|T_i+S|\nonumber \\
&\ge& |T|+|W|\frac{|S|}{2}.\label{X11}
\end{eqnarray*}
It follows that $W=\{j\}$, form some $j$. Since $T+S$ is aperiodic, we have by Kneser's Theorem, $|T_j+S|\ge |T_j|+|S|-1$.
Since $$|T|+|S|-1=|T+S|\ge \sum _{i\neq j}|T_i+S|+|T_j|+|S|-1,$$ the result follows \end{proof}


\begin{lemma}\label{transfer}
Let $S$ and $T$ be finite subsets of an abelian group $G$ such that
 $0\in S\cap T$ and $S+T$ is aperiodic. Let $H$ be a finite subgroup and let $\phi:G\mapsto G/H$ denotes the canonical morphism.
Let   $S=\bigcup \limits_{0\le i\le u}S_i$ be a $H$--decomposition such that $\phi(S_0), \cdots ,\phi(S_u) $ is a progression with difference $d$  and  that $S\setminus S_u$
 is $H$--periodic. Then there is a decomposition $T=\bigcup \limits_{0\le i\le t}T_i$, such that  $\phi(T_0), \cdots ,\phi(T_t) $ is a progression with  difference $d$ and
 such that $T\setminus T_t$
 is $H$--periodic.
\end{lemma}

The proof is an easy exercise.



\begin{lemma}\label{tpowers}
Let $S$ and $T$ be subsets of a finite  abelian group $G$, generated by $S$ such that
 $S+T$ is aperiodic,  $0\in S\cap T$ and $|S+T|=|S|+|T|-1$. Then
 $T^S-S$ is aperiodic and  $|T^S-S|=|T^S|+|S|-1$.

 \end{lemma}

 \begin{proof} The set $T^S-S$ is aperiodic by Lemma \ref{lee}. Clearly $T^S-S\subset G\setminus T.$ Thus $|T^S-S|\le |G|-|T|=|T^S|+|S|-1$. By Kneser's Theorem we have
$|T^S-S|\ge |T^S|+|S|-1.$\end{proof}

\section{The $\frac{2|G|}3$--Theorem }

The following result encodes efficiently the critical pair Theory.

\begin{theorem}\label{twothird}
Let $S$ be  a finite  generating subset  of
 an abelian group $G$ of order $\ge 2$.
Let $T$ be  a finite  subset  of
 $G$ such that
 $|S|\le |T|$,  $S+T$ is aperiodic,  $0\in S\cap T$ and
 $$ \frac{2|G|+2}3\ge |S+T|= |S|+|T|-1.$$

 Let $H$ be a hyper-atom   of $S$  and let $\phi :G\mapsto G/H$ denotes
the canonical morphism.  Then
\begin{itemize}
\item $S$ and $T$ are $H$--quasi-periodic,
\item  $\phi(S)$ and $\phi(T)$ are arithmetic progressions with the same
 difference.
 \end{itemize}

\end{theorem}

\begin{proof}

Set $|G|=n$, $h=|H|$,
$|\phi(S)|=u+1$, $|\phi(T)|=t+1$ and
$q=\frac{n}{h}$.

We have $|S|\le \frac{|S|+|T|}2\le \frac{n+1}3<\frac{n+1}2$.

Assume first that  $|H|=1$. By Corollary \ref{vosper}, $S$ is an arithmetic progression.
The result holds clearly, observing that $S$ generates $G$.

Assume now that $|H|\ge 2.$
Take  $H$--decompositions $S=\bigcup \limits_{0\le i\le u}S_i$,
$T=\bigcup \limits_{0\le i\le t}T_i$ and $S+T=\bigcup _{0\le i \le k}E_i$.

Without loss of generality, we shall assume that
\begin{itemize}
  \item $0\in S_0$ and  $|S_0|\ge \cdots \ge |S_u|$.
  \item $|E_{t+1}|\ge \cdots \ge |E_k|$ and $T_i+S_0\subset E_i$, for all $0\le i \le t$.
\end{itemize}

For $0\le i \le t$, we put $\alpha _i=|E_i|-|T_i|$.
By the definition we have
$u|H|=|H+S|-|H|=\kappa _2(S)\le |S|-1.$ It follows that for all $u\ge
j\geq 0$
\begin{equation}\label{plein}
|S_{u-j}|+\cdots +|S_u|\ge j|H|+1
\end{equation}

Put $P=\{i \in [0,t] :  |E_i|=h\},$
   $W=\{i \in [0,t] :  |E_i|<h\}.$

By (\ref{plein}), $|S_0|>\frac{h}2$.
Thus $\subgp{S_0}=H$, by Lemma \ref{prehistorical}.
We have by (\ref{olson}), $|T_i+S_0|\ge |T_i|+\kappa _1(S_0) \ge |T_i|+\frac{|S_0|}{2},$ for all $i\in W$. Therefore
\begin{eqnarray}|T+S|
&\ge& \sum _{i\in W}|T_i+S_0|+\sum _{i\in P}|E_i|+\sum _{t+1\le i \le k}|E_i|
\nonumber \\&\ge&
 |T|+|W|\frac{|S_0|}{2}+\sum _{i\in P}\alpha _i+\sum _{t+1\le i \le k}|E_i|.\label{X11}
\end{eqnarray}

{\bf Claim 0}  $ t+1+u\le q.$

Suppose the contrary. By Lemma \ref{prehistorical}, $\phi(S\setminus S_u)+\phi(T)=G/H$. In particular $k+1=q$ and
 $|E_i|\ge |S_{u-1}|$, for all $i$. Therefore we have by 
 (\ref{plein}),
\begin{eqnarray*}
|S+T|
&\ge &(t+1) |S_0|+(q-t-1)|S_{u-1}|\\
&=& (2t+2-q)|S_0|+(q-t-1)(|S_0|+|S_{u-1}|)\\
&=& (2t+2-q)\frac{uh+1}{u+1}+2\frac{uh+1}{u+1}(q-t-1)\\
&\ge&\frac{un+3}{u+1},
\end{eqnarray*}
noticing that $q< t+1+u\le 2t+1.$ Therefore $u=1$ and $q=t+1$.

Since  $\frac{n}3>|S|-1\ge \kappa _2(S)\ge h=\frac{n}q$, we must have $q\ge 4.$

We must have $|W|\ge 3, $  since otherwise by (\ref{X11}), $$|T+S|\ge (q-2)h+2|S_0|\ge(q-1)h+1=n-h+1>
\frac{2n}3+1, $$ a contradiction. We must have $|W|=3,$ since otherwise by (\ref{X11}), $|T+S|\ge |T|+2|S_0|\ge |T|+|S|, $ a contradiction. Then $|P|=q-|W|\ge 4-3=1$.

Since $\subgp{S}=G$ and $u=1$, we have
$\subgp{\phi(S)}=\subgp{\phi(S_1)}=G/H$, and hence there a $\gamma \in P$ such that $T_\gamma+S_1\subset E_j,$ for some $j\in W$.

By Lemma \ref{prehistorical}, $ |T_\gamma| +|S_j|\le h$. In particular, $\alpha _\gamma \ge |S_1|$. By (\ref{X11}), $|T+S|\ge |T|+|S_1|+\frac{3|S_0|}2>|T|+|S|,$
a contradiction.
Hence $ t+1+u\le q.$

{\bf  Claim} 1. $|\phi (S+T)|=|\phi (S)|+|\phi (T)|-1.$

By Claim 0, (\ref{cosetgraph}) and (\ref{eqisoper0}), we have
$$k+1=|\phi(S+T)|\ge \min(q,t+1+u)=t+1+u.$$

 Also, we have $(t+1)h\ge |T|\ge |S|> \kappa _2(S)=uh$, and hence
$t\ge u.$  By Lemma \ref{cosetgraph},
$\kappa _1(\phi (S))=  |\phi (S)|-1$.

By Proposition \ref{strongip}, there
a  subset $J\subset
[0,t]$ with  $|J|=u$ and a family $\{ mi ;i\in J\}$ of
integers in $[1,u]$
  such that
   $T+S $ contains the $H$--decomposition $(\bigcup \limits_{0\le i\le
  t}T_i+S_0)\cup (\bigcup \limits_{ i\in
  J}T_i+S_{mi})$.

We have $k= t+u,$ since otherwise
\begin{eqnarray*}
|S+T| &\ge& \sum _{i\in J}|T_i+S_0|+\sum _{i\notin  J}|T_i+S_0|+\sum _{i\in J}|T_i+S_{mi}|+\sum _{i\ge t+u+1}|E_i|\label{AP1}\\
&\ge& u|S_0|+|T|+|S_{u}|\ge |T|+|S|, \label{AP2}
\end{eqnarray*}
a contradiction. Thus $$|\phi (S+T)|=|\phi (S)|+|\phi (T)|-1.$$

{\bf Claim 2} $|E_{k-1}|\ge |S_{u-1}|$.

By Theorem \ref{hyperatom}, $\phi(S)$ is  an arithmetic progression or a Vosper subset. Let us show that
\begin{equation}
|\phi (T)+\phi (S\setminus S_0)|\ge t+u \label{omit}
\end{equation}

Notice that (\ref{omit}) is obvious if $\phi(S)$ is  an arithmetic progression
and follows  by Lemma \ref{vominus} if $\phi(S)$  is a Vosper subset. Claim 2 follows now.

{\bf Claim 3} If $u\ge 2$ then $q-1\ge t+u+2$.

Assume   $u\ge 2$.  We must have
 \begin{equation}  |P|\ge 2.\label{EQP>1}\end{equation}
Suppose the contrary.
By Lemma \ref{prehistorical},  $|T_i|<\frac{h}{3}<\frac{|S_0|}2$ for every $i\in W$. We have using Claim 2 and (\ref{plein}), \begin{eqnarray*}
2|T|>|S+T|&\ge& \sum _{i\in W}|T_i+S_0|+ \sum _{i\in P}|E_i|+ \sum _{t+1\le i\le k}|E_i|\nonumber\\
&\ge& \sum _{i\in W}2|T_i|+ |P||H|+|S_{u-1}|+ |S_u|\label{EQV2T0}
\\ &\ge& \sum _{i\in W}2|T_i|+ 2h+1>2|T|,\label{EQV2T}
\end{eqnarray*}
a contradiction.
Take  $U\subset P$ with $|U|=2.$
We have using (\ref{plein}),
 $$ 2|S_0|\ge |S_0|+|S_{u-1}|\ge \frac{2}{3}(|S_u|+|S_{u-1}|+|S_{u-2})\ge \frac{4h+2}3.$$
 By (\ref{omit}), we have $|E_i|\ge |S_{u-1}|$ for all $i\le k-1$. Clearly   $|E_i|\ge |S_0|$, for all $0\le i \le t$.
The claim must hold since otherwise we have using Claim 2 and (\ref{plein}),
  \begin{eqnarray*}
|S+T|&\ge&  \sum _{i \in U} |E_i|+(t-1)|S_0|+(u-1)|S_{u-1}|+|S_u|\\
&=& 2h+ (t-2)|S_0|+(u-2)|S_{u-1}|+(|S_0|+|S_{u-1}|+|S_u|)\\
&=&  2h+(t-u)|S_0|+(u-2)(|S_0|+|S_{u-1}|)+(|S_0|+|S_{u-1}|+|S_u|)\\
&\ge& 2h+ (t-u)\frac{2h+1}3+\frac{(4h+2)(u-2)}3+2h+1\\
&>&(t+u+2)\frac{2h}3+1\ge q\frac{2h}3+1=\frac{2n}{3}+1,
\end{eqnarray*}
 a contradiction.


By Claim 1 and Claim 3, $\phi(S)$ can not be a Vosper subset. By Theorem \ref{hyperatom}, $\phi(S)$ is  an arithmetic progression.
By Lemma \ref{AP1} and since $t+1+u<q$ for $u\ge 2$,
 $\phi(T)$ is an arithmetic progression with the same difference as $\phi (S)$.

Now we shall reorder the $S_i$'s And $T_i$'s using the modular progression structure.

Take  $H$--decompositions $S=\bigcup \limits_{0\le i\le u}B_i$, $ T=\bigcup \limits_{0\le i\le t}A_i$
and a $H$--decomposition $S+T=\bigcup _{0\le i \le k}R_i$.
Since $-d$ is a difference of $\phi(S)$,  we may assume $0\in B_0$ and that
\begin{enumerate}
\item  $\phi(B_0), \cdots , \phi(B_u)$ is an arithmetic progression with difference $d$ and $|B_0|\ge |B_u|$.
\item  $\phi(A_0), \cdots , \phi(A_t)$ is an arithmetic progression with difference $d$.
\item $A_i+B_0\subset R_i$, for all $0\le i \le t$;
\item $A_{t}+B_{i}\subset R_{t+i}$, for all $1\le i \le u$.
\end{enumerate}
We shall put $Y=\{i \in [0,t] :  |R_i|<h\}.$
Since $|B_0|\ge |B_u|$, we have using (\ref{plein}), $|B_0|>\frac{h}2$.
Thus $\subgp{B_0}=H$, by Lemma \ref{prehistorical}.
We have by (\ref{olson}), $|A_i+B_0|\ge |A_i|+\kappa _1(B_0) \ge |A_i|+\frac{|B_0|}{2},$ for all $i\in Y$.
Thus

\begin{eqnarray}|T+S|&\ge & \sum _{0\le i \le t}|A_i+B_0|+\sum _{1\le i \le u}|A_t+B_i|\\
&\ge & |T|+|Y|\frac{|B_0|}{2}+\sum _{1\le i \le u}|A_t+B_i|.\label{Y11}
\end{eqnarray}

By (\ref{Y11}),
we have
$|T+S|\ge|T|+|Y|\frac{|B_0|}{2}+|S\setminus B_0|,$
and hence $|Y|\le 1.$

{\bf  Claim} 4. $Y=\emptyset$. Suppose the contrary. Then $Y=\{r\}$, for some $0\le r \le t$. Assume first $r<t$. By Lemma \ref{prehistorical},
$h\ge |A_r|+|B_0|$. Thus
\begin{eqnarray*}
|S+T|&\ge&|R_r|+  th+ |A_t+(B\setminus B_{0})|\\
&\ge & |B_{0}|+|A_r|+|B_{0}|+(t-1)h+|A_t|+ \sum _{1\le i\le u-1}|B_{i}|\\
&\ge& |T|+|S|-|B_u|+|B_0|\ge |S|+|T|,
\end{eqnarray*} a contradiction.
Then $r=t$. By Lemma \ref{prehistorical},
$h\ge |A_t|+|B_0|$. Also $|R_t|\ge |A_{t-1}|$. Hence
\begin{eqnarray*}
|S+T|&\ge& th+|R_t|+  |A_t+(B\setminus B_{0})|\\
&\ge & |A_t|+|B_{0}|+(t-1)h+|A_{t-1}|+ \sum _{1\le i\le u}|B_{i}|\\
&\ge& |T|+|S|,
\end{eqnarray*} a contradiction.

Let us show that  $|E_i|=h,$ for all $i\le k-1$.

Suppose  that  there is an  $r\le k-1$
with $|E_r|<h$. By Claim 4, $t+1\le r$.
 Since $(A_{t}+B_{r})\cup (A_{t-1}+B_{r+1})\subset R_r$, we have using (\ref{plein}),
  $2h\ge |A_t|+|B_r|+|A_{t-1}|+|B_{r+1}|\ge |A_t|+|A_{t-1}|+h+1$, by Lemma \ref{prehistorical}.
  Thus   $|T+H|-|T|\ge 2h-(|A_t|+|A_{t-1}|)\ge h+1.$
  Now $|S+T|\ge |T+H|+|S\setminus B_0|\ge |T|+h+1+|S|-|B_0|>|T|+|S|,$ a contradiction.

Since $S+T$ is aperiodic, the set  $B_u+A_t$ is aperiodic. By Kneser's Theorem
$|B_u+A_t|\ge |B_u|+|A_t|-1.$ Now have
\begin{eqnarray*}
|S|+|T|-1=|S+T|&\ge& (t+u)h+|A_t+B_{u}|\\
&\ge & th+|B_{u}|+th+|A_t|-1+ \sum _{1\le i\le u-1}|B_{i}|\\
&\ge& |T|+|S|-1
\end{eqnarray*}
Thus $|S\setminus B_u|=uh$ and $|T\setminus A_t|=th$.\end{proof}

Notice that  the subgroup in Theorem \ref{twothird} depends only one of the sets (namely $S$),
while the subgroup in Kemperman Structure Theorem depends on $S$ and $T$.

\section{The  Structure Theory}

A pair $\{A,B\}$ of subsets of an abelian group $G$ will be called a  {\em weak pair} if one of
the following conditions holds:
\begin{itemize}
\item [(WP1)] $\min(|A|,|B|)=1$;
\item [(WP2)] There is $d\in G$, with order $\ge |A|+|B|-1$, such that $A$ and $B$ are arithmetic progressions with difference
$d$;

\item [(WP3)] $A$ is aperiodic and there is a finite subgroup  $H$ and $g\in G$ such that $A,B$ are contained in some
$H$--cosets and   $g-B=H\setminus A$;
\item [(WP4)] There is a subgroup  $H$ such that $A,B$ are contained in some $H$-cosets and $|A|+|B|=|H|+1
$, and  moreover there is  $c\in G$ such that
$|(c-A)\cap B|=1.$

\end{itemize}

An elementary pair is a  pair $\{A,B\}$ satisfying one of the conditions (WP1), (WP2), (SP3) and  (SP4), where

\begin{itemize}
  \item $SP3$ = "$WP3$ and  for every $c\in G$, such that
$|(c-A)\cap B|\neq 1,$
  \item  $SP4$ = "$WP4$ and  for every $d\in G\setminus \{c\}$,
$|(d-A)\cap B|\neq 1.$
\end{itemize}


The notion of a weak  pair was
suggested by Kemperman \cite{kempacta} (end of section 5, Page 82) in order to formulate an easier description.
The reader could use the following Lemma, implicit in Kemperman's work \cite{kempacta}, if he wants to work with elementary pairs.

\begin{theirlemma} (Kemperman)
Let $A,B$ be  a
weak pair.

  There is a nonzero subgroup $H$   and  $H$-quasi-periodic decompositions
 $A=A_0\cup A_1$ and $B=B_0\cup B_1$ such that $(A_1,B_1)$ is an
elementary pair and $|(\phi(a_1+b_1)-\phi(A))\cap \phi (B)|=1$,
 where $\phi : G\mapsto G/H$  is the canonical morphism and
 $a_1\in A_1$ and $b_1\in B_1$.
\label{strongpair}
 \end{theirlemma}

 The proof of the above lemma is implicit Kemperman's work \cite{kempacta}.
It can be done within one page or less using Kneser's Theorem. The reader may
also refer to the Appendix of Lev's paper \cite{levkemp}.

\begin{theorem}\label{finalcor}

Let $T$ and $S$ be  a finite  subset  of an abelian
 $G$ such that $\subgp{S\cup T}=G$.
 Assume that $|S|\le |T|$,   $0\in S\cap T$ and
 $$ |S+T|= |S|+|T|-1.$$
Also assume that $S+T$ is aperiodic and that  $\{S,T\}$ is not a weak pair.
 \begin{itemize}
   \item[(i)] If $G\neq \subgp{S}$ then  $T$ is $\subgp{S}$--quasi-periodic;
   \item[(ii)] If $G= \subgp{S}$ and $G\neq  \subgp{U}$ then  $S$ and $T$ are $U$--quasi-periodic, where $U=T^S-a$, for some
   $a\in T^S$;
    \item[(iii)] If $G= \subgp{S}= \subgp{U}$ then
  there is a
 proper subgroup $H$ of  $L$ such that $S$ and $T$ are $H$--quasi-periodic modular progressions.
 Moreover $H$ is  either a hyper-atom of $S$ or a hyper-atom of $U$.

 \end{itemize}

\end{theorem}

\begin{proof}
Since $\{S,T\}$ is not a weak pair, we have $|S|\ge 2$.
Put $L=\subgp{S}$. If $L\neq G$, then $T$ is $L$--quasi-periodic by Lemma \ref{nongenerating} and (i) holds. So we may assume without loss of generality
that $S$ generates $G$.
We have $|T^S|\ge 2$, otherwise putting $G=(T+S)\cup \{v\}$, we have
clearly $v-S\subset G\setminus T$. Actually we have equality by the relation $|T+S|=|T|+|S|-1$.
 Then $\{T,S\}$ is a weak pair, a contradiction.

Notice that $|S|+|T|+|T^S|=|S|+|T|+(|G|-|T+S|)=|G|+1$.


\begin{itemize}
  \item If $|T|\le |T^S|$, then $|S|+|T|\le \frac{2(|S|+|T|+|T^S|)}3=\frac{2|G|+2}3$, and (iii) holds by Theorem \ref{twothird}.
  \item  $|T|> |T^S|$,  then $|S|+|T^S|\le \frac{2(|S|+|T|+|T^S|)}3=\frac{2|G|+2}3$.

   Assume first $|U|\ge |S|$.
By Theorem \ref{twothird}, $S$ a quasi-periodic modular $H$--progression, where $H$ is the hyper-atom of $S$. By Lemma \ref{transfer},
 $T$ a quasi-periodic modular $H$--progression. Thus (iii) holds.

 Assume now $|U|<|S|$ and that $U$ generates $G$.
 By Theorem \ref{twothird}, $S$ a quasi-periodic modular $H$--progression, where $H$ is the hyper-atom of $U$. By Lemma \ref{transfer},
 $T$ a quasi-periodic modular $H$--progression. Thus (iii) holds.

\item If  $\subgp{U}\neq G$, then
 by Lemma \ref{nongenerating}, $S$ and $T$ are $\subgp{U}$--quasi-periodic. Thus (ii) holds.

\end{itemize}

 The proof is complete. \end{proof}


Following a suggestion of Kemperman \cite{kempacta}, we formulate the main classical critical pair result in  the following way:

\begin{corollary}\label{KST} (Kemperman Structure Theorem \cite{kempacta})

Let $A,B$ be finite subsets of an abelian group $G$ with $|G|\ge 2.$

Then  the following conditions are equivalent:

\begin{itemize}
\item [(I)] $|A+B|= |A|+|B|-1,$ and moreover $|(c-A)\cap B|=1$ for some
$c$ if $A+B$ is periodic.

\item [(II)]  There is a nonzero subgroup $H$   and  $H$-quasi-periodic decompositions
 $A=A_0\cup \cdots \cup A_u$ and $B=B_0\cup \cdots \cup B_t$ such that $(A_0,B_0)$ is an
 elementary pair and $|(\phi(A_0+B_0)-\phi(A))\cap \phi (B)|=1$,
 where $\phi : G\mapsto G/H$  is the canonical morphism.
 \end{itemize}
\end{corollary}

\begin{proof}

The implication (II) $\Rightarrow $ (I) is quite easy.
Suppose that (I) holds.
Without loss of  generality we may take $|A|\ge |B|$,  $0\in A\cap B$ and $G=\subgp{A\cup B}$.
 Let $Q$ denotes the period of  $A+B$ and let $\psi : G\mapsto G/Q$  denotes the canonical morphism.

{\bf Case} 1. $|Q|=1.$

Clearly  (II) holds with $H=\{0\}$
 if $\{A,B\}$ is a weak pair. Suppose that $\{A,B\}$ is not a weak pair.
 By Theorem \ref{finalcor},  there is a proper subgroup $H$ such that $A$ and $B$ are    $H$--quasi-periodic.
 Take a minimal such a group. By $\rho _H(X)$ we shall mean the unique non-full $H$--coset's
trace on $X$ if such a coset exists.

 The pair $\{\rho(A), \rho(B)\}$ must be a
weak pair, by  Theorem \ref{finalcor}. Since $A+B$ is aperiodic, we have $|(\phi(\rho(A)+\rho(B))-\phi(A))\cap \phi (B)|=1$,
 where $\phi : G\mapsto G/H$  is the canonical morphism.

{\bf Case} 2. $|Q|\ge 2.$

Take two $Q$--decompositions $A=A_0\cup \cdots \cup A_u$ and $B=B_0\cup \cdots \cup B_t$ such that
$c\in A_0+B_0$.
Observe that $\psi (c)=\psi (A_0)+\psi (B_0)$ has a unique
expression. Hence by Scherck's Theorem \ref{scherk}, $|\psi(A)+\psi
(B)|\ge |\psi(A)|+|\psi (B)|-1=t+u+1.$ We must have $|\psi(A)+\psi (B)|=
t+u+1,$ since otherwise
$|A+B|=|\psi(A+B)||H|\ge |A|+|B|.$

 By
Lemma \ref{prehistorical}, we have $$|A_0|+|B_0|\le |H|+1.$$
Now we have
\begin{eqnarray*}
\sum _{0\le i \le u}|A_i|+\sum _{0\le i \le t}|B_i|-1&=& |A+B|\\&=&|H|(u+t+1)\ge |H|(u+t)+|A_0|+|B_0|-1\ge
|A|+|B|-1,
\end{eqnarray*}

observing
that $|A|+|B|-1 = |A+B| = |\psi(A+B)||H|=|A+H|+|B+H|-|H|.$  It follows that

\begin{itemize}
  \item
$|A_1|=\cdots =|A_u|=|B_1|=\cdots =|B_t|=|H|,$
  \item $|A_0|+|B_0|\le |H|+1.$
\end{itemize}
Therefore  $\{A_0,B_0\}$ is a weak  pair.
\end{proof}

We shall now prove Lev's  Structure Theorem.

\begin{corollary}\label{LEVK} (Lev's Structure Theorem)\cite{levkemp}

Let A and B be finite nonempty subsets of a non-zero abelian group
$G$, with  $|A|\le |B|$,  and  $|A+B|=|A|+|B|-1$. Suppose that either $A+B$ is aperiodic or there is a $c$ with $|(c-A)\cap B| = 1$.

Then  there exists a finite  proper subgroup $H\subset G$ such that $A$ and $B$ are $H$--quasi-periodic and $\{ \phi(A), \phi (B)\}$ is a weak pair.

\end{corollary}

\begin{proof}

Without loss of generality we may assume that $0\in A\cap B$ and $\subgp{A\cup B}=G$.
Put $L=\subgp{B}$.

Let $P$ be the  period of
$A+B$  and let $\sigma : G\mapsto G/P$  be the
canonical morphism.

{\bf Case} 1. $|P|=1$.
\begin{itemize}
  \item [(i)] If  $G\neq L$ then  $A$ is $L$--quasi-periodic and $|\phi(B)|=1$, where $\phi : G\mapsto G/L$
is the canonical morphism, by Theorem \ref{finalcor}.
  \item [(ii)]If $G= L \neq \subgp{U}$, then  by Theorem
\ref{finalcor},  $B$ and  $A$ and $\subgp{U}$--quasi-periodic, where $U=A^B-a,$ for some $a\in A^B.$ Also  by Theorem
\ref{finalcor}, $\phi(B+A)=G/\subgp{U}$, where $\phi : G\mapsto G/\subgp{U}$
is the canonical morphism.
  \item [(iii)]If $G= L= \subgp{U}$, then
  there is a
 proper subgroup $H$  such that $B$ and $A$ are $H$--quasi-periodic modular progressions, by Theorem
\ref{finalcor}.
\end{itemize}

In order to show that $\{\phi(B),\phi (A)\}$ is a weak pair, it would be enough to observe that there is a uniquely
representable element in $\phi(B)+\phi (A)$. This follows easily since $A+B$ is aperiodic and since $A$ and $B$ are quasi-periodic.

{\bf Case} 2. $A+B$ is periodic.

By Case 1, there is subgroup $K$ of $G/P$ such that $\sigma (A)$ and $\sigma (B)$
are $K$--quasi-periodic and
$\{\tau (\sigma (A)),\tau (\sigma (B))\}$
is a weak pair,  where $\tau : G/K\mapsto G/P/K$. Put $\psi=\tau\circ \sigma$ and $Q=\psi ^{-1}(K)$.
Clearly $A$ and $B$ are $Q$--quasi-periodic. Clearly
$\{\psi(A),\psi (B)\}$ is a weak pair.\end{proof}

{\bf Acknowledgement}. The author wishes an anonymous referee for  valuable comments on a preliminary draft related to this work.

\end{document}